\documentstyle{article}
\newtheorem{theorem}{Theorem}
\newcommand{\hio}{ H^{\infty} (\Omega)}
\newcommand{\whot}{H^1 (T,B)}
\newcommand{\hit}{ H^{\infty} (T)}
\newcommand{\Dl}{\Delta}
\begin{document}
\title{Jean Bourgain's analytic partition of unity \\ via holomorphic 
martingales}
\author{Paul F. X. M\"uller\thanks{Supported by FFWF Pr.Nr. JP 90061} \\
Division of Mathematics, California 
Institute of Technology\\ CA 91125, Pasadena.\\and
 \\Institut f\"ur Mathematik J. Kepler Universit\"at \\ Linz,Austria.}  
\date{24 September 1992}
\maketitle 
\begin{abstract}
Using 
stopping time arguments on holomorphic martingales we present a soft way
of constructing J. Bourgain's analytic partitions of unity. Applications to
Marcinkiewicz interploation in weighted Hardy spaces are discussed.
\end{abstract}
\section{Introduction}
In his 1984 Acta Mathematica paper Jean Bourgain
derives new Banach space properties of $H^{\infty}$ and the disc 
algebra from the existence of the following analytic partition of unity:
\begin{theorem}[J. Bourgain]
 Given $f$, a strictly positive integrable function on ${\bf T}$
with $\int f(t) dt  = 1$ and $ 0<\delta<1$ then, there exist functions $\tau_j, 
\gamma_j \in H^{\infty}(T) $ and positive numbers $c_i$ such that:
\begin{enumerate}   
\item $ \| \gamma_j \|_{\infty} < C$
\item $ \sum  |w_j|<C$
\item $|\tau_j|f <c_j$
\item $\sum c_j \|\tau_j\|_1 <\delta^{-C}$ 
\item $\int |1-\sum \gamma_j \tau_j^2 |fdt<\delta .$
\end{enumerate}
\end{theorem}
Here I whish to present a soft way to this 
construction 
which results from the use of {\bf probabilistic tools} such as 
holomorphic martingales.
I should like to point out here that a proof for the existence
of analytic partitions of unity --much simpler than J. Bourgain's --
 has been given recenty by Serguei Kislyakov. See [K1] and [K2].\par
In [K3] S. Kislyakov derived J. Bourgain's result on p-summing operators
from the following weighted Marcinkiewicz decomposition. 
\begin{theorem}[S. Kislyakov]
For any positive weight b on $\bf{ T}$ there exists a weight $B \geq b$ and
$\int Bdt < C \int bdt$ so that for any $\lambda >0$ and $f \in \whot $
there exists $g \in \hit $ and $ h \in \whot $ satisfying.
\begin{enumerate}
\item $f=g+h$
\item $||g||_{\infty} \le \lambda $
\item $\int hBdt \le C \int _{\{|f|>\lambda\}} |f|Bdt.$
\end{enumerate}
\end{theorem} 
Up to small perturbations we shall obtain 
a {\bf stochastic version} of  Kislyakov's decomposition 
which allowes us to prove  
the following: 
\begin{theorem}[J. Bourgain]    
For any 2-summing operator $S$ on the disc algebra  
and any $2<q <\infty$ the q-summing norm satisfies the
interpolation inequality
$$ \pi_q (S) \le C_q \pi_2 (S)^{2 \over q}||S||^{1- {2 \over q}}.$$
\end{theorem}
A very elegant proof of this interpolation inequality has been given 
by Gilles Pisier who used vectorvalued $H^1$ spaces. 
See [P].  
\section{The main result}
Holomorphic martingales were introduced by N. Varopoulos in [V]. 
They are stable under stopping times, and generalize analytic
functions on the unit circle. This connection has lead to probabilistic 
proofs of several results in Analysis, including Carleson's corona
thoerem [V], the existence of a logmodolar Banach algebra having no
analytic structure [C] and P.W. Jones's interpolation theorems
between $H^1$ and $H^\infty$ [M1,M2].\par
This paper is not selfcontained! We freely use notations and definitions from 
[V] without further explanation.
\begin{theorem}
Given $\Delta$, a strictly positive integrable function on $(\Omega, P)$
with $\int \Delta d P = 1$ and $ 0<\delta<1$ , there exist functions $w_j, 
\theta_j \in \hio$ and positive numbers $c_i$ such that:
\begin{enumerate}   
\item $ \| \theta_j \|_{\infty} < C$
\item $ \sum  |w_j|<C$
\item $|w_j|\Delta <c_j$
\item $\sum c_j \|w_j\|_1 <\delta^{-C}$ 
\item $\int |1-\sum \theta_j w_j^2 |\Delta dP<\delta $
\end{enumerate}
\end{theorem}
Probability offers a {\bf soft way of constructing}  the functions $\theta_j$ 
so that the verification if (5) becomes much easier than
in J. Bourgain's proof. See  [B, pp. 11, 12]. The probabilistic concept
will be merged with analytic tools, such as Havin's lemma, which we use in 
the following form, due to Bourgain:      
\begin{theorem} For every measurable subset $E$ of $\Omega$ and $0<\epsilon 
<1 $ 
there exist functions $ \alpha , \beta \in  \hio $ such that:
\begin{enumerate}
\item $|\alpha|+|\beta| \le 1$
\item $|\alpha-\frac{1}{5}| < \epsilon$ on $E$.  
\item $|\beta | < \epsilon $ on $E$.
\item $\| \alpha \| _1 < C| \log \epsilon |^2 P(E)$
\item $\| 1-\beta \| _2 < | \log \epsilon | P(E)^{1 \over 2}.$
\end{enumerate} 
\end{theorem} 
{\bf Proof of Theorem 4:} We shall first determine a new weight:  
Let d be the outer function,
so that $|d|=\Dl$ and put
$$ A(\Dl):=   sup_t |E(d|{\cal F}_t)| $$
then we let $$ \Delta _1 = \sum_{n=0}^{\infty} A^n (\Dl) (C2)^{-n} $$    
where $ C$ is determined by Varopoulos' inequality: For $d \in H^1(\Omega)$
$$ \int sup_t |E(d|{\cal F}_t)| < C \|d\|_1 $$
Clearly, this construction gives,
      
\begin{enumerate}   
      \item $ A(\Delta _1) < \Delta _1 3C $
\item $ \Dl < C\Delta _1 $
\item $ \int \Delta_1 dP<C \int \Dl dP. $
      \end{enumerate}
We next define holomorphic partitions of unity: Let $\Psi$
be the outer function so that $|\Psi| =\Delta_1$.  
Consider now the stopping times $ \tau _0 =0$ and 
$$ \tau_j:=inf\{t> \tau_{j-1} : |E(\Psi | {\cal F}_{t})| > M^j \} $$
to define $\Psi_i:=E(\Psi| {\cal F}_{\tau_{j}})$ and $ d_j :=\Psi_{j+1} -
\Psi_{j} 
,$ elements of $\hio$ for or which, obviously the identity
$$ 1=\frac{E(\Psi)}{\Psi}+ \sum_{j=0}^{\infty} \frac{d_j}{\Psi}$$
holds. The summands of the above expression will be our choice of $\theta_i$:
Indeed we define
 $ \theta_{-1} := \frac{E(\Psi)}{\Psi}$ and 
           $ \theta_i := \frac{d_{i}}{\Psi}$  for $ i=0,1,2,\dots$
Obviously we obtain $$\|\theta_i\|_{\infty} \le C. $$ 
 Havin's lemma allows us to truncate the above partition of unity:
We apply it to sets $E_i := \{ \Psi{^*} > M^i\}$ and denote 
the resulting functions by $\alpha_i, \beta_i$. Then define for $i=-1,0,1,...$
$$ w_i := 5\alpha_i \prod_{s=8}^{\infty} \beta_{i+s}^{s}$$  
{\bf Verification of property 
(5).}  \hfill \break
We first eliminate the weight $\Delta$: 
$$ \int |1-\sum_{i=-1}^{\infty} \theta_i w_i^2 |\Delta dP=\int |\sum_{i=-
1}^{\infty} 
\theta_i(1-w_i^2)|\Delta 
dP \le 
 2\sum_{i=-1}^{\infty} \int |d_i(1-w_i)|dP $$
Using the inequality
$$ |1-\prod z_i| \le \sum |1-z_i| $$
which holds for complex numbers in the closed unit disc, 
we get the following upper bound for the obove sum of integrals:
$$\sum_{i} \int |\ d_i (1-5\alpha_i)
       + \sum_{s>8} \int |d_i (1-\beta_{i+s}^{s})| dP
$$ 
The martingale differences $ d_i$ are supported on $E_i$ and bounded by 
$M^{i+1}$. Therefore we obtain a domination by:
   $$\sum_i \int_{E_i} |(1-5\alpha_i)M^{i+1}
       + \sum_{s>8} M^{i+1} \int_{E_i} M^{i+1} |(1-\beta_{i+s}^{s})| dP
$$ 
Invoking the estimates from Havin's Lemma and applying  Cauchy-Schwarz'
inequality give the following estimates:
   $$\epsilon M C       +\log(\epsilon^{-1}) 
\sum_{s>8} \sum_{i=-1}^{\infty} M^{i+1}sP(E_i)^{1 
\over 2} P(E_{i+s})^{1 \over 2}  
$$ 
Again by Cauchy-Schwarz we dominate the above sum by:
   $$\epsilon M C       +\log(\epsilon^{-1})  \sum_{s>8} M^{1-\frac{s}{2}}s
(\sum_ {i=-1}^{\infty} M^{i}P(E_i))^{1 \over 2 } ( \sum_{i=-1}^{\infty} M^{i+s}
P(E_{i+s})) ^{1 \over 2} \, \le 
$$ 
   $$\epsilon M C       +\log(\epsilon^{-1})  \sum_{s>8} M^{1-\frac{s}{2}}sC
$$ 
This is what we want if $\epsilon$ is chosen of order $M^{-2}$ and 
$M:=\delta^{-1}.$ \hfill\break
Havin's lemma, repeatedly applied,  
gives the following statements:
\begin{enumerate}
\item $ \sum  |w_i|<C$
\item $|w_i|\Delta <M^{i+8}$
\item $\sum M^{i}  \|w_i\|_1 \le\sum M^{i}  \| \alpha_i\|_1 \le\sum M^{i} 
P(E_i)|\log \epsilon | $
\end{enumerate}
To finish the proof it is now enough to take $c_i=M^i$ 
\section{Reduction of J. Bourgain's partition of unity}
To obtain Bourgain's original result, we lift the density $f$ from
$T$ to $\Omega$ construct a new weight together with holomorphic 
partitions of unity there and project the solutions back to $\Omega.$
This is done by norm-one operators $$M:H^p(T) \to H^p(\Omega)$$  and
$$N:H^p(\Omega) \to H^p(T)$$ so that $Id=NM$, and $N(M(f)F)=fN(F)$. 
\hfill\break  
{\bf Proof of theorem 1:}
Apply Theorem 4 to the density $\Delta:=Mf. $ Let $$g_i := N(\theta_i w_i^2)$$
We define  $\gamma_i$ to be the inner factor of $g_i$ and put
$$ \tau_i := a_i^{1 \over 2}$$
where $a_i$ denotes the outer factor of $g_i.$
Using our main result it is easy 
to verify conditions $1) \dots...5)$ of Bourgain's theorem.

\section {Truncating functions in weighted $H^p$}
Here we combine stopping times and holomorphic partitions of unity to obtain 
Marcinkiewicz decomposition in weighted Hardy spaces. \hfill \break 
Although the next theorem looks terribly complicated, it   
simply states 
that up to a (reasonable) change of density and up to a small error,
interpolation is 
possible in weighted Hardy spaces. 
\begin{theorem} For any density $\Delta$ on $\Omega$ 
and $ \delta > 0 $ there exists $\phi \in \hio $ so that for any 
\begin{enumerate}
\item $\| \phi \|_{\infty} < C$
\item  $\Delta_1 > \Dl$ and $\int \Delta_1 dP < \delta^{-C} \int \Delta dP $
\item $\int |1-\phi| \Delta dP < \delta \| \Delta \|_1$
\item For any $f \in  H^q(\Omega, \Delta_1) $ and any $\lambda > 0 $
there exists $g \in \hio $ and $h \in H^2 (\Omega, \Delta_1)$
satisfying 
\begin{enumerate}
\item $  f\phi = g+h $
\item $\|g\|_{\infty} \le \lambda$
\item $\int |h|^2 \Delta_1dP < C_q \lambda^{2-q} \int |f^q| \Delta_1 dP $

\end{enumerate}
\end{enumerate}
\end{theorem}
Proof: Let  $w_i \in \hio $ and $\theta_i \in \hio $ be given by Theorem 4 .
Then we define:$$
\phi := \sum \theta_i w_i^2$$
$$ \Delta _1 := \Delta + \sum c_i |w_i|  $$
$$ f_i = w_if$$
Now we use the stopping time 
$$ \tau_j:=inf\{t : |E(f_j | {\cal F}_{t})| > \lambda \} $$
to define $g_j:=E(f_j| {\cal F}_{\tau_{j}})$ and $h_j:=f_j -g_j$.
By the stability property of holomorphic martingales these functions
are certainly holomorphic and satisfy       
\begin{enumerate}
\item $||g_j||_{\infty} \le \lambda $
\item $\int h_jdP \le 2\int _{\{|f_j ^*|>\lambda\}} |f_j|dP$
\end{enumerate}
Now, using partitions of unity we glue these partial solutions together
$$g:=\sum g_j w_j \theta_j $$ and  
$$ h:= \sum h_jw_j \theta_j $$
Then clearly $$g+h=\sum(g_j+h_j)v_j=f \phi $$
and $$\|g\|_{\infty}< sup_{j}\|g_j\|_{\infty} \|\sum |w_j| \|_{\infty}<
\lambda C $$                   
The estimate for $\int |h|^2 \Delta_1 dP$ follows a well established 
pattern,
which has been carefully presented in the central chapter of Wojtaszczyk's book.
See [W, Ch III.I]. \hfill \break Property 3) of theorem 4 implies that:
$$\int |h|^ 2 \Delta_1 dP \le C \sum \int |h_j|^2 |w_j|^2  \Delta_1 dP$$
The last sum can be estimated, using the interplay between the partitions 
and the density, by                                        
$$ \sum c_j \int|h_j|^2dP \le \sum c_j \int _{\{|f_j ^*|>\lambda\}} |f_j|^2dP$$
Using H\"older's inequality for conjugate indices $r,s$ we estimate the above 
expression by 
$$(\sum c_j\|f_j\|_{2r}^{2r})^{1 \over r}
(\sum \lambda^{-2r}c_j\|f_j^*\|_{2r}^{2r})^{1 \over s} $$

As the martingale maximal function is bounded in $L^{2r}(\Omega,P)$ we obtain 
an upper bound proportional to
$$(\int |f|^{2r}\sum c_j|w_j|dP)^{1 \over r}
(\sum \lambda^{-2r}c_j\|f_j\|_{2r}^{2r})^{1 \over s} $$
Specializing $r={q \over 2}$ this product is finally dominated by a 
constant,depending on $q$  
times
$$  \lambda^{2-q} \int |f|^q\Delta_1 dP. $$
\section{Reduction of J.Bourgain's interpolation inequality}
For the 2- summing operator $S$ there exists a positive probability
measure on {\bf T} so that
$$ \|Sx\| \le \pi_2(S) (\int |x| d\mu)^{1 \over 2}  \,\, for \, \, x \in A$$
Without loss or generality we may assume that $\mu $ is absolutely
continuous w.r.t. Lebesgue measure, i.e.,
$$ d\mu = fdt. $$ 
Consquently for $b \in \hio $ the operator $U =SN$ satisfies
$$ \|Ub\|\le \pi_2 (S) (\int |b| \Delta dP)^{1 \over 2} $$    
whrere $ \Delta =Mf$.\hfill \break
{\bf Proof of Theorem 3:} Let $ 0<\delta <1 $ be given. Theorem 6 applied to
the density $\Delta$ shows that $U$ can be split into $ U = U_1 + R_1$
so that $$ \pi_q(U_1) \le \delta ^{-C}\pi_2 (S)^{2 \over q}   \|S\|^{1-{2 \over
q}}$$                                                                         
and $$ \pi_2 (R_1) \le \delta \pi_2(S).$$
where $U_1 g=U(b\phi)$ and $R_1b=U(b(1-\phi))$. Indeed fix $b \in H^q(\Omega, 
\Delta_1)$
of norm one in that space. Then according to Theorem 6 for $\lambda =\|S\|  
^{q \over 2}\pi_2 (S)^{-{2 \over q}}  $ we find a 
Marcinkiewicz decomposition of $\phi b$ into 
$$ \phi b = g+h.$$                         
Therefore $$ \|Ub\phi\| \le \|Ug\| + \|Uh\|< \|S\|\lambda + \pi_2 
(S)\lambda^{1-{ q\over 2}} 
 \le \|S\|  ^{1 -{2\over q}} \pi_2 (S)^{{2 \over q}}. $$
As for the error term we have
$$ \|U(b(1-\phi))\| \le \pi_2(S)( \int |b|^2 |1-\phi|^2 \Delta
dP)^{1 \over 2}  $$
Using property 3), 4) and 5) of theorem 4 gives
$$ \int |1-\phi|^2  \Delta dP  \le C \int |1-\phi| \Delta dP \le C \delta \int 
\Delta dP$$             
and 
$$   \int  \Delta_1 dP \le \delta^{-C}   \int \Delta dP$$ 
We therefore obtained the correct estimates for $U_1$ and $R_1.$ To
finish the proof 
of theorem 3, we now iterate the above decompositionand observe that
$S=UM$.
\section{References} 
\begin{enumerate}
\item  
[B1] J. Bourgain, New Banach space properties of the disc algebra and 
$H^\infty$,Acta. Math. 152 (1984) 
\item [B2] \dots           , Bilinear forms on $H^\infty$ and bounded
bianalytic 
functions, Trans.Amer. Math. Soc 286 (1984)
\item [C] K.Carne, The algebra of bounded holomorphic martingales, 
J.F.Aanal. 45 (1982)  
\item   [K1] S. Kislyakov, Absolutely summing operators on the disc algebra,
St.Petersburg Math.J 3  (1991)
\item [K2] \dots , Extensions of (p,q) summing operator on the disc algebra
with an appendix on Bourgain's analytic partition, Preprint 1990.
\item [K3]   \dots,          Truncating functions in weighted $H^p$ and 

two theorems of
J. Bourgain, Preprint Uppsala University, 1989.
\item  [M1] P.F.X.M\"uller, Holomorphic martingales and interpolation between
Hardy spaces, J. d'Analyse Math. Jerusalem, to appear.
\item [M2] \dots ,   Holomorphic martingales and interpolation
between Hardy spaces: The 
complex method,preprint 1992.
\item [P] G. Pisier, A simple proof of a theorem of J. Bourgain, preprint 1990
\item [V] N. Varopoulos, Helson Szego Theorem $A_p$ functions for Brownian
motion and several variables, J.F.Anal. 39 (1980)\hfill
\item [W] P. Wojtaszczyk, Banach spaces for analysts, Cambridge
Univ. Press 1991\hfill
\end{enumerate}

\end{document}